\documentclass[10pt]{article} \usepackage{amssymb,amscd}
\usepackage[all]{xy} \usepackage[cp1255]{inputenc}
\usepackage[hebrew,english]{babel} \usepackage{amsmath}
\usepackage{xspace}

\makeindex
\author{Oz Ben-Shimol
 \\ \\
\small{Department of Mathematics,
Bar-Ilan University,} \\
\small{Ramat-Gan, Israel}}
\date{\small{\today}}
\title{Centers associated with the Borel subalgebra of the \\ general linear Lie algebra\thanks{Work is supported by ISF grant \#170/12.}}
\begin{document}
\maketitle
\newcommand{\gr}{\operatorname{gr}}
\newcommand{\si}{\operatorname{si}}
\newcommand{\ad}{\operatorname{ad}}
\newcommand{\tr}{\operatorname{trace}}
\newcommand{\depth}{\operatorname{depth}}
\newcommand{\Mod}{\operatorname{mod}}
\newcommand{\Ht}{\operatorname{ht}}
\newcommand{\Cl}{\operatorname{Cl}}
\newcommand{\Spec}{\operatorname{Spec}}
\newcommand{\Sl}{\operatorname{sl}}
\newcommand{\Span}{\operatorname{Span}}
\newcommand{\Char}{\operatorname{char}}
\newcommand{\Id}{\operatorname{Id}}
\newcommand{\Der}{\operatorname{Der}}
\newcommand{\cpf}{\operatorname{pf_{c}}}
\newcommand{\codim}{\operatorname{codim}}
\newcommand{\rest}[1]{\lvert_{#1}}

\newcommand{\ZZ}{{\mathbb Z}} \newcommand{\RR}{{\mathbb R}}
\newcommand{\NN}{{\mathbb N}} \newcommand{\QQ}{{\mathbb Q}}
\newcommand{\fg}{{\mathfrak{g}}} \newcommand{\HH}{{\mathbb H}}
\newcommand{\fn}{{\mathfrak{n}}}
\newcommand{\fhn}{{\mathfrak{\hat{n}}}} 
\newcommand{\fr}{{\mathfrak{r}}}
\newcommand{\fm}{{\mathfrak{M}}}
\newcommand{\DD}{{\mathfrak{D}}}
\newcommand{\KK}{{\mathfrak{K}}}
\newcommand{\BB}{{\mathfrak{b}}}
\newcommand{\FF}{{\mathbb F}}
\newcommand{\cL}{{\mathcal{L}}}

\newcommand{\eee}{\hfill$\Box$}

\thispagestyle{empty}
\baselineskip=0.55cm
\noindent
\thispagestyle{empty}
\begin{center}
    \large \textbf{Abstract}
\end{center}

We consider a Borel subalgebra $\fg$ of the general linear algebra and its subalgebra $\BB$ which is a Borel subalgebra of the special linear algebra, over arbitrary field. 
Let $\cL\in\{\fg, \BB\}$.
We establish here explicit realizations of the center $Z(\cL)$ and semi-center $Sz(\cL)$ of the enveloping algebra, the Poisson center $S(\cL)^{\cL}$ and Poisson semi-center $S(\cL)^{\cL}_{\si}$ of the symmetric algebra. We describe their structure as commutative rings and establish isomorphisms $Z(\cL)\cong S(\cL)^{\cL}$,  $Sz(\cL)\cong S(\cL)^{\cL}_{\si}$
\subsection*{Introduction}
\ \ \ \ Let $\fg$ be the Lie algebra of the upper triangular matrices with $n$ rows and $n$ columns with elements in an \emph{arbitrary} field $F$ of characteristic $p\geq 0$. 
$\fg$ is a Borel subalgebra of the general linear Lie algebra.
For $p$ which does not divide $n$, we consider its Lie subalgebra $\BB$ which is a Borel subalgebra of the special linear algebra.
For $\cL\in\{\fg,\BB\}$, let $S(\cL)$ be its symmetric algebra, $U(\cL)$ its enveloping algebra. 
The Poisson center of $S(\cL)$ is denoted by $S(\cL)^{\cL}$ (the set of invariants of the $\cL$-module $S(\cL)$ with respect to the adjoint representation), 
and the Poisson semi-center of $S(\cL)$ is denoted by $S(\cL)^{\cL}_{\si}$ (4.1).  
The center of $U(\cL)$ is denoted by $Z(\cL)$ and the semi-center of $U(\cL)$ is denoted by $Sz(\cL)$ (5.1). 
\\ \indent
In Section 1 we assume $p>0$ and establish an explicit realization of $S(\fg)^{\fg}$ by generators and relations.
It is shown that $S(\fg)^{\fg}$ is a Cohen-Macaulay ring and it is a complete intersection ring (1.11) for $p=2$. 
This realization leads to an explicit realization of $S(\BB)^{\BB}$ by generators and that $S(\BB)^{\BB}$ is a Cohen-Macaulay ring.
\\ \indent
In Section 2, still with $p>0$, we use the realizations in Section 1 and the fact that $S(\fg)$ is the graded algebra of $U(\fg)$, to obtain the following results: 
an explicit realization of $Z(\fg)$, an explicit isomorphism of $S(\fg)^{\fg}$ onto $Z(\fg)$ (thus $Z(\fg)$ is a Cohen-Macaulay ring and it is a complete intersection ring for $p=2$),  
realization of $Z(\BB)$ by generators, and an explicit isomorphism of $S(\BB)^{\BB}$ onto $Z(\BB)$ (thereby $Z(\BB)$ is a Cohen-Macaulay ring).
\\ \indent
In Section 3, we apply reduction modulo $p>0$ in conjunction with the realization of $S(\fg)$ in Section 1, 
to obtain a realization of $S(\fg)^{\fg}$, $Z(\fg)$,  $S(\BB)^{\BB}$ and $Z(\BB)$ for $p=0$ (that is, over fields of zero characteristic). 
It is shown that $S(\fg)^{\fg}$ and $Z(\fg)$ are isomorphic as polynomial rings of one variable. 
Furthermore, we prove $S(\BB)^{\BB}=Z(\BB)=F$. Loosely speaking, $U(\BB)$ has no center.
\\ \indent
In Sections 4,5, we get the analogue results for the semi-centers. 
We get realizations of $S(\fg)^{\fg}_{\si}$, $Sz(\fg)$, $S(\BB)^{\BB}_{\si}$, $Sz(\BB)$ and isomorphisms $S(\fg)^{\fg}_{\si}\cong Sz(\fg)$, $S(\BB)^{\BB}_{\si}\cong Sz(\BB)$.
It follows that for $p>0$ the algebras $S(\fg)^{\fg}_{\si}$,  $S(\BB)^{\BB}_{\si}$ (hence $Sz(\fg)$, $Sz(\BB)$) are complete intersection rings, and if $p=0$ then $S(\fg)^{\fg}_{\si}$ (hence $Sz(\fg)$) is a polynomial ring of $n$ variables, while  $S(\BB)^{\BB}_{\si}$ (hence $Sz(\BB)$) is a polynomial ring of $n-1$ variables. 
\\ \indent
We end this introduction by mentioning related known results.
\\ \indent
Over algebraically closed fields of characteristic zero, there exists an isomorphism of the Poisson center onto the center (the so called \emph{Duflo isomorphism}) [4, Theorem 10.4.5].
There exists, as well, an isomorphism of the Poisson semi-center onto the semi-center [8].  
\\ \indent
A generating set for $Z(\fg)$ is introduced, without proof, by V.V.Panyukov [7].
\\ \indent
Over the field of complex numbers, using completely different methods, A. Josef [5] shows that $Sz(\BB)$ is a polynomial ring of $n-1$ variables.
\\ \indent
A description of $Sz(\BB)$ for $p>0$ is given by G. Vernik [9];  he describes it as a polynomial ring over the center of a certain subalgebra of $U(\BB)$.  
\newpage
\textbf{Notation}
$$
\begin{array}{lll}
	F & \text{-} & \text{the ground field}. \\ 
	p & \text{-} & \text{the characteristic of} \ F \ (p\geq 0). \\ 
	e_{s,t} & \text{-} &  \text{the} \ n\times n \ \text{matrix} \ (\alpha_{k,l}) \ \text{such that} \ \alpha_{k,l}=0 \ \text{for} \ (k,l)\neq(i,j) \ \text{and} \ \alpha_{i,j}=1. \\ 
	\fm & \text{-} &  \text{the matrix} \ (e_{i,j})_{1\leq i,j\leq n}. \\ 
	\fg & \text{-} & \text{the Lie algebra of the upper triangular} \ n\times n \ \text{matrices over} \ F. \\
           &             & \fg=\Span_{F}\{e_{i,j} \ | \ 1\leq i\leq j\leq n\}. \\ 
	\fn & \text{-} & \text{the Lie algebra of the strictly upper triangular} \ n\times n \ \text{matrices over} \ F.  \\
          &               & \fn=\Span_{F}\{e_{i,j} \ | \ 1\leq i< j\leq n\}. \\
	\{\fn\} &         & \text{the underlying set of} \ \fn. \\
	\BB   & \text{-} & \text{a Borel subalgebra of} \ \Sl_{n}(F). \\ 
         &               & \BB=\Span_{F}(\{\fn\}\cup\{e_{i,i}-e_{n,n} \ | \ 1\leq i\leq n-1\}), \ \\ 
	\cL & \text{-} & \text{a Lie algebra}. \\
	S(\cL) & \text{-} & \text{the symmetric algebra of} \ \cL. \\ 
	S_{p}(\fg),S_p(\BB) & \text{:} & \text{1.2}, \ \text{1.13}. \\ 
	S(\cL)^{\cL} &\text{-} & \text{the Poisson center of} \ S(\cL). \\ 
	S(\fg)^{\fn} &\text{=} & \{f\in S(\fg) \ | \ \ad{x}(f)=0, \ x\in\fn\} \\ 
	S(\cL)^{\cL}_{\si} &\text{-} & \text{the Poisson semi-center of} \ S(\cL).\ (\text{4.1}) \\ 
	U(\cL) & \text{-} & \text{the enveloping algebra of} \ \cL. \\ 
	Z(\cL) & \text{-} & \text{the center of} \ U(\cL). \\ 
	Z_{p}(\fg),Z_{p}(\BB) & \text{:} & \text{(2.1, 2.2), 2.8}. \\
	Sz(\cL) & \text{-} & \text{the semi-center of} \ U(\cL). \ \text{(5.1)} \\ 
	Q(A) & \text{-} & \text{the quotient field of a domain A}. \\ 
	h   & = & \lfloor\frac{n-1}{2}\rfloor. \\ 
	C(k) & \text{-} & \text{the determinant of the} \ k\text{-th right upper block of}\ \fm, \ \ 1\leq k\leq\lfloor\frac{n}{2}\rfloor. \\ 
     D(k) & = & \displaystyle{\sum_{i=1}^{k}e_{i,i}+e_{n-i+1,n-i+1}}, \ \  1\leq k\leq h. \\ 
     T_{k}(i,j)  & \text{-} & \text{the determinant of the matrix which obtains from the} \ k\text{-th right upper} \\ 
               &          & \text{block of} \ \frak{M} \ \text{by replacing its} \ i \text{-th row with the last} \ k\text{-th tuple of the} \\
               &          &  j \text{-th row of} \ \fm  \ \ (1\leq i\leq k \ \text{and} \ k+1\leq j\leq n-k). \\ 
     S_{k}(i,j) & \text{-} & \text{the determinant of the matrix which obtains from the} \ k\text{-th right upper} \\ 
               &          & \text{block of} \ \frak{M} \ \text{by replacing its} \ i \text{-th column with the first} \ k\text{-th tuple of the} \\
               &          &  (n-j+1) \text{-th row of} \ \fm  \ \ (1\leq i\leq k \ \text{and} \ k\leq j\leq n-k). \\ 
	T(k)  & = & \displaystyle{\sum_{i=1}^{k}\sum_{j=k+1}^{n-k}e_{i,j}T_{k}(i,j)}, \ \  1\leq k\leq h. \\ 
     M(k) & = & C(k)D(k)+T(k), \ \  1\leq k\leq h, \ \ M(\frac{n}{2})=0 \ \text{for an even n}. \\ 
     M_{\BB}(k) & \text{:} & \text{1.13.} \\ 
	c_{0},z_{0}	& \text{:} & \text{1.2, 2.1.}\\
	c(k,l),z(k,l)          & \text{:} &  \text{1.2, 2.1.} \\
	c_{\BB}(k,l),z_{\BB}(k,l)  & \text{:} &  \text{1.13, 2.8.}\\
	d(r) & \text{-} &  \text{the set consisting of the elements of the} \ r\text{-th upper diagonal of} \ \fm. \\
            & &  \text{That is,} \ d(r)=\{e_{i,j} \ | \ j-i=r\}.	 \\ 
	\hat{d}(r) & \text{-} & \text{the subset of} \ d(r) \ \text{consisting of the elements which do not belong to the} \\ 
                  & & \text{anti-diagonal of} \ \fm. \ \text{That is,} \ \hat{d}(r)=d(r)\setminus\{e_{s,n-s+1} \ | \ s=1,\ldots,n\}.
\end{array}
$$
\newpage
\subsection*{1. \ The Poisson center of $\mathbf{S(\fg)}$}
\textbf{1.1.} \ In this section, $\Char(F)=p>0$.
\vskip 0.1cm\noindent
\textbf{1.2.} \textbf{Generators.} 
$$
\begin{array}{lll}
	S_{p}(\fg) & = & F[e_{i,j}^{p} \ | \ 1\leq i\leq j\leq n]. \\ \\ 
     c_{0}        & = & \displaystyle{\tr(\frak{M})=\sum_{i=1}^{n}e_{i,i}}. \\ 
                      &    &  \text{This element generates the center of} \ \fg. \\ \\
     c(k,l)         & = & C(k)^{p-l}M(k)^{l}, \ \  1\leq k\leq h, \ 0\leq l\leq p-1. \\
                      &    & \text{In particular} \  c(k,0)=C(k)^{p}\in S_{p}(\fg).
\end{array}
$$
\textbf{1.3.} \ In the sequel we shall need the following basic Lemma.
\vskip 0.1cm\noindent
\textbf{Lemma.} \
Let $X=(x_{i,j})$ be an $r\times r$ matrix over a commutative ring $A$.
If $D\in\Der(A)$ then $D(\det(A))$ can be computed column by column.
That is, if $X=(X_1,\ldots,X_r)$ where $X_1,\ldots,X_r$ are the columns of $X$, then
$$
    D(\det(X))=\sum_{i=1}^{r}\det(X_1,\ldots,X_{i-1},D(X_{i}),X_{i+1},\ldots,X_r)
$$
where $D(X_i)=(D(x_{1,i}),D(x_{2,i}),\ldots,D(x_{r,i}))^{T}$.
\vskip 0.1cm\noindent
\textbf{1.4. \ Proposition.} \\
\textbf{a.} \ $C(k),M(k)\in S(\fg)^{\fn}$. \\ 
\textbf{b.} \ $c(k,l)\in S(\fg)^{\fg}$.
\vskip 0.1cm
\textbf{Proof.} \ The fact $C(k)\in S(\fg)^{\fn}$ is proved in [1]. 
Let separate the underlying set of $\fm$ to six subsets as follows:
$$
\begin{array}{ll}
	\fm_1=\{e_{s,t} \ | \ 1\leq s<t\leq k \} & \ \fm_4=\{e_{s,t} \ | \ k+1\leq s<t\leq n-k \} \\ \\ 
     \fm_2=\{e_{s,t} \ | \ 1\leq s\leq k, \ k+1\leq t\leq n-k \} & \ \fm_5=\{e_{s,t} \ | \ k+1\leq s\leq n-k<t\leq n \} \\ \\ 
     \fm_3=\{e_{s,t} \ | \ 1\leq s\leq k, \ n-k+1\leq t\leq n \} & \ \fm_6=\{e_{s,t} \ | \ n-k+1\leq s<t\leq n \}
\end{array}
$$
We have 
$$
\ad{e_{s,t}}(T_{k}(i,j))=\left\{\begin{array}{ll}
		-\delta_{si}T_{k}(t,j) & e_{s,t}\in\fm_1 \\ \\
           \delta_{si}\delta_{tj}C(k) & e_{s,t}\in\fm_2 \\ \\
		\delta_{tj}T_{k}(i,s) & e_{s,t}\in\fm_4 \\ \\ 
		0                              & e_{s,t}\in\fm_3\cup\fm_5\cup\fm_6
\end{array}\right.
$$
Since
$$
	\ad{e_{s,t}}(T(k))=\sum_{i=1}^{k}\sum_{j=k+1}^{n-k}\left[(\delta_{ti}e_{s,j}-\delta_{sj}e_{i,t})T_{k}(i,j)+e_{i,j}\ad{e_{s,t}}(T_{k}(i,j))\right],
$$
we get
$$
	\ad{e_{s,t}}(T(k))=\left\{\begin{array}{ll}
		0 			               &  e_{s,t}\in\fm_1\cup\fm_3\cup\fm_4\cup\fm_6 \\ \\
           e_{s,t}C(k)             &  e_{s,t}\in\fm_2 \\ \\
           -\sum_{i=1}^{k}e_{i,t}T_{k}(i,s) & e_{s,t}\in\fm_5
\end{array}\right..
$$
It is very easy to see that
$$
	\ad{e_{s,t}}(D(k))=\left\{\begin{array}{ll}
		0 			               &  e_{s,t}\in\fm_1\cup\fm_3\cup\fm_4\cup\fm_6 \\ \\
           -e_{s,t}             &  e_{s,t}\in\fm_2 \\ \\
           e_{s,t}               & e_{s,t}\in\fm_5
\end{array}\right..
$$
Since
$$
	\ad{e_{s,t}}(M(k))=C(k)\ad{e_{s,t}}(D(k))+\ad{e_{s,t}}(T(k)),
$$
we get 
$$
	\ad{e_{s,t}}(M(k))=0 \ \ \ \text{if} \ \ \ e_{s,t}\in\fm_1\cup\fm_2\cup\fm_3\cup\fm_4\cup\fm_6
$$
and
$$
	\ad{e_{s,t}}(M(k))=e_{s,t}C(k)-\sum_{i=1}^{k}e_{i,t}T_{k}(i,s)   \ \ \ \text{if} \ \ \ e_{s,t}\in\fm_5.
$$
Therefore, to complete the proof of \textbf{a.} we have to show that 
$$
	\sum_{i=1}^{k}e_{i,t}T_{k}(i,s)=e_{s,t}C(k)   \ \ \ \text{if} \ \ \ e_{s,t}\in\fm_5.
$$
Let $M(i,r)$ be the $(i,r)-$th minor of the matrix of $C(k)$, $1\leq i,r\leq k$. 
Using the expansion of $T_{k}(i,s)$ by the $i$-th row yields
$$
	\sum_{i=1}^{k}e_{i,t}T_{k}(i,s)=\sum_{i=1}^{k}e_{i,t}\sum_{r=1}^{k}(-1)^{i+r}e_{s,n-k+r}M(i,r)	
$$  
$$
	=\sum_{r=1}^{k}e_{s,n-k+r}\sum_{i=1}^{k}(-1)^{i+r}e_{i,t}M(i,r)
$$
$$
	=\sum_{r\neq t-(n-k)}e_{s,n-k+r}\sum_{i=1}^{k}(-1)^{i+r}e_{i,t}M(i,r)
$$
$$
    +e_{s,t}\sum_{i=1}^{k}(-1)^{i+(t-(n-k))}e_{i,t}M(i,t-(n-k)).
$$
Note that the expression $\displaystyle{\sum_{i=1}^{k}(-1)^{i+r}e_{i,t}M(i,r)}$ is the expansion of the determinant by the $r-$th column, 
of the matrix which obtains from the matrix of $C(k)$ by replacing its $r-$th column with its $(n-k+t)-$th column. 
Therefore, if $r\neq t-(n-k)$, this expression is a determinant of a matrix with two equal columns and therefore zero. 
Also, the expression $\displaystyle{\sum_{i=1}^{k}(-1)^{i+(t-(n-k))}e_{i,t}M(i,t-(n-k))}$ is the expansion of $C(k)$ by its $(t-(n-k))-$th column.
\newline\indent
To prove \textbf{b.}, a description of the action of $\ad{e_{s,s}}$ is required.   
It is easy to verify that
$$
\ad{e_{s,s}}(T_{k}(i,j))=\left\{\begin{array}{ll}
		(1-\delta_{si})T_{k}(i,j) & 1\leq s\leq k \\ \\
           \delta_{sj}T_{k}(i,j) & k+1\leq s\leq n-k \\ \\
		-T_{k}(i,j) & n-k+1\leq s\leq n 
\end{array}\right..
$$
Since
$$
	\ad{e_{s,s}}(T(k))=\sum_{i=1}^{k}\sum_{j=k+1}^{n-k}\left[(\delta_{si}e_{s,j}-\delta_{sj}e_{i,s})T_{k}(i,j)+e_{i,j}\ad{e_{s,s}}(T_{k}(i,j))\right],
$$
we get
$$
	\ad{e_{s,s}}(T(k))=\left\{\begin{array}{ll}
		T(k)   &  1\leq s\leq k \\ \\
           0       &  k+1\leq s\leq n-k \\ \\
           -T(k)  & n-k+1\leq s\leq n
\end{array}\right..
$$
Obviously, $\ad{e_{s,s}}(D(k))=0$, and
$$
	\ad{e_{s,s}}(C(k))=\left\{\begin{array}{ll}
		C(k)   &  1\leq s\leq k \\ \\
           0       &  k+1\leq s\leq n-k \\ \\
           -C(k)  & n-k+1\leq s\leq n
\end{array}\right..
$$
Since
$$
	\ad{e_{s,s}}(M(k))=\ad{e_{s,s}}(C(k))D(k)+\ad{e_{s,s}}(T(k)),
$$
we get
$$
	\ad{e_{s,s}}(M(k))=\left\{\begin{array}{ll}
		M(k)   &  1\leq s\leq k \\ \\
           0       &  k+1\leq s\leq n-k \\ \\
           -M(k)  & n-k+1\leq s\leq n
\end{array}\right..
$$
Finally, 
$$
	\begin{array}{lll}
		\ad{e_{s,s}}(c(k,l)) & = &\ad{e_{s,s}}(C(k)^{p-l}M(k)^{l}) \\ \\
                                       & = & (p-l)C(k)^{p-l-1}\ad{e_{s,s}}(C(k))M(k)^{l} \\ \\
                                       & + & lC(k)^{p-l}M(k)^{l-1}\ad{e_{s,s}}(M(k))=0.
\end{array}
$$
\eee
\vskip 0.1cm\noindent
\textbf{1.5.} \textbf{Relations.} \ For $0\leq i,j\leq p-1$ let $r(i,j)$, $s(i,j)$ be the unique integers such that $0\leq s(i,j)\leq 1$, $0\leq r(i,j)\leq p-1$ and $i+j=ps(i,j)+r(i,j)$. 
Then for $k=1,\ldots,h$ we have
\begin{equation}
	c(k,i)c(k,j)=c(k,r(i,j))C(k)^{p(1-s(i,j))}M(k)^{ps(i,j)}.
\end{equation}
\vskip 0.1cm\noindent
\textbf{1.7. Theorem.} \ $\mathit{Q(S(\fg)^{\fg})=Q(S_{p}(\fg)[c_{0},c(k,1) \ | \ k=1,\ldots,h])}$. 
\vskip 0.1cm
\textbf{Proof.} \ \textit{Step 1. \ The field extension} 
$$
	\mathit{Q(S_{p}(\fg)[c_{0},c(k,1) \ | \ k=1\ldots,h])/Q(S_p(\fg))}
$$ 
\textit{is of degree} $\mathit{p^{h+1}}$: 
\\ \indent 
$c_0\not\in S_{p}(\fg)$, $c_{0}^{p}\in S_p(\fg)$, thus $Q(S_p(\fg)[c_0])/Q(S_p(\fg))$ is of degree $p$. \\
Let $i,j,l,k$ be integers, $1\leq k\leq h$, $1\leq l<k$, $1\leq i\leq l$, $l+1\leq j\leq n-l$, and $\DD=\ad{e_{n-k+1,k}}$.
Then $\DD{D(l)}=\DD{C(l)}=\DD{c_0}=0$. Also, $\DD{T_{l}(i,j)}=\delta_{j,k}T_{l}(i,n-k+1)$, hence $\DD{T(l)}=0$.
It follows that $\DD{M(l)}=0$ and therefore 
\begin{equation}
	\DD{S_{p}(\fg)[c_0,c(1,1),\ldots,c(l,1)]}=0.
\end{equation}
Now, $\DD{D(k)}=0$, and $\DD{C(k)}\in S(\fg)$ for if $e_{s,t}$ is an entry of the underlying matrix of $C(k)$ then $1\leq s\leq k$, $n-k+1\leq t\leq n$ hence 
$$
	\frak{D}e_{s,t}=\delta_{k,s}e_{n-k+1,t}-\delta_{n-k+1,t}e_{s,n-k+1}\in S(\fg).
$$
We shall see that $\DD{T(k)}\not\in S(\fg)$. Let $e_{s,t}$ be an entry of the underlying matrix of $T_k(i,j)$ where $1\leq i\leq k$, $k+1\leq j\leq n-k$. 
Then $\DD{e_{s,t}}=\delta_{k,s}e_{n-k+1,t}-\delta_{n-k+1,t}e_{s,k}$.
Therefore, $\DD{e_{s,t}}\not\in S(\fg)$ if and only if $s=j$ and $t=n-k+1$, that is, if and only if $e_{s,t}$ is the $i$-th entry of the first column of the underlying matrix of $T_k(i,j)$.
Apply $\DD$ to $T_{k}(i,j)$ column by column (Lemma 1.3) and conclude that $\DD{T_{k}(i,j)}=e_{k,j}\alpha(i,j)+\beta(i,j)\not\in S(\fg)$ for some $\alpha(i,j),\beta(i,j)\in S(\fg)$.
Now, $\DD{e_{i,j}}=\delta_{k,i}e_{n-k+1,j}$ hence
$$
	\DD{T(k)}=\sum_{j+1}^{n-k}e_{n-k+1,j}T_{k}(k,j)+\sum_{i=1}^{k}\sum_{j=k+1}^{n-k}e_{i,j}[(e_{j,k}\alpha(i,j)+\beta(i,j)],
$$
which is certainly not an element of $S(\fg)$.
\\ \noindent
So we have $\DD{M(k)}=D(k)\DD{C(k)}+\DD{T(k)}\not\in S(\fg)$ and therefore 
$$
	\DD{c(k,1)}=-C(k)^{p-2}\DD{C(k)}M(k)+C(k)^{p-1}\DD{M(k)}\neq 0.
$$
From $(2)$ it follows that $c(k,1)\not\in S_{p}(\fg)[c_0,c(1,1),\ldots,c(k-1,1)]$. \\
Clearly $c(k,1)^{p}\in S_{p}(\fg)[c_0,c(1,1),\ldots,c(k-1,1)]$, hence the field extension
$$    
	Q(S_{p}(\fg)[c_0,c(1,1),\ldots,c(k,1)])/Q(S_{p}(\fg)[c_0,c(1,1),\ldots,c(k-1,1)])
$$
is of degree $p$ for every $k=1,\ldots,h$. 
\\ \indent
\textit{Step 2.} We proceed as in [1, Proposition 1.13] when we are dealing here with fields instead of division algebras.
\\ \indent
Let $r,i,j,k$ be integers such that $1\leq r\leq n-1$, $1\leq k\leq\lfloor\frac{n}{2}\rfloor$, $1\leq i\leq k$, $k\leq j\leq n-k$. 
\\
For each $e_{i,j}\in d(r)$ we have \\
(1) \ If $e_{i,j}$ is lying above the anti-diagonal of $\fm$ then 
$$
\begin{array}{l}
\ad{T_{i}(i,j)}(e_{i,j})\neq 0 \ \text{and} \\ \\ 
\ad{T_{i}(i,j)}(d(r)\setminus\{e_{i,j}\}\cup d(r+1)\cup\ldots\cup d(n-1))=0.
\end{array} 
$$
(2) \  If $e_{i,j}$ is an element of the anti-diagonal of $\fm$ or lying below it, then 
$$
\begin{array}{l}
\ad{S_{i}(1,n-j+1)}(e_{i,j})\neq 0 \  \text{and} \\ \\
\ad{S_{i}(1,n-j+1)}(d(r)\setminus\{e_{i,j}\}\cup d(r+1)\cup\ldots\cup d(n-1))=0. 
\end{array}
$$
Denote by $\{\fn\}$ the underlying set of $\fn$. 
So $\{\fn\}=\bigcup_{r=1}^{n-1}d(r)$.
It follows that the field extension $Q(S(\fg)^{\fg}[\{\fn\}])/Q(S(\fg)^{\fg})$ is of degree $p^{\deg{\fn}}=p^{\frac{n(n-1)}{2}}$. 
\\ \indent
Now, for $k=1,\ldots,\lfloor\frac{n}{2}\rfloor$ we have $\ad{C(k)}(e_{k,k})=-C(k)$, $\ad{C(k)}(e_{l,l})=0$ when $l<k$, and by Proposition 1.2.a. $\ad{C(k)}(S(\fg)^{\fg}[\{\fn\}])=0$.
It follows that the field extension 
$$
	Q(S(\fg)^{\fg}[\{\fn\}\cup\{e_{k,k} \ | \ k=1,\ldots,\lfloor\frac{n}{2}\rfloor\}])/Q(S(\fg)^{\fg}[\{\fn\}])
$$ is of degree $p^{\lfloor\frac{n}{2}\rfloor}$.
Therefore, the degree of 
$$
	Q(S(\fg)^{\fg}[\{\fn\}\cup\{e_{k,k} \ | \ k=1,\ldots,\lfloor\frac{n}{2}\rfloor\}])/Q(S_p(\fg))
$$ is at least 
$$
	p^{h+1}p^{\frac{n(n-1)}{2}}p^{\lfloor\frac{n}{2}\rfloor}=p^{\frac{n(n+1)}{2}}=p^{\dim{\fg}},
$$
which is the degree of $Q(S(\fg))/Q(S_{p}(\fg))$. 
By degree considerations we must have 
$$
	Q(S(\fg)^{\fg})=Q(S_{p}(\fg)[c_{0},c(k,1) \ | \ k=1,\ldots,h]).
$$  
\eee
\vskip 0.1cm\noindent  
\textbf{1.8. Lemma.} \ $\mathit{S_{p}(\fg)[c_0]}$ \textit{is a polynomial ring.} 
\vskip 0.1cm
\textbf{Proof.} \ By changing basis in $\fg$ we get that 
$$
	S_{p}(\fg)=F[c_{0}^p,e_{i,j}^p \ | \ 1\leq i\leq j\leq n, (i,j)\neq(n,n)].
$$
Hence 
$$
	S_p(\fg)[c_0]=F[c_0,e_{i,j}^p \ | \ 1\leq i\leq j\leq n, (i,j)\neq(n,n)].
$$ 
$c_0$ is integral over $S_{p}(\fg)$, thus the Krull dimension of the affine ring $S_p(\fg)[c_0]$ is equal to the amount of its generators. 
$S_p(\fg)[c_0]$ is a polynomial ring.
\eee
\vskip 0.1cm\noindent 
\textbf{1.9.} \ Let $t_0,t(k,l)$,\ $k=1,\ldots,h$, $l=1,\ldots,p-1$, be algebraically independent elements over  $S_{p}(\fg)$.
Set $t(k,0)=C(k)^{p}$.
Consider the polynomial ring
$$
	R=S_{p}(g)[t_0,t(k,l) \ | \ k=1,\ldots h, \ l=1,\ldots p-1] 
$$
and the following polynomials of $R$ (which are corresponding to relations (1)): 
$$
\begin{array}{l}
	f_0=t_{0}^{p}-c_{0}^{p} \\ \\
	f_{k}(i,j)=t(k,i)t(k,j)-C(k)^{p(1-s(i,j))}M(k)^{ps(i,j)}t(k,r(i,j)), \\ \\
     \ \ \ \ \ \ 1\leq i\leq j\leq p-1, \ k=1,\ldots,h. 
\end{array}
$$
Let $I$ be the ideal of $R$ generated by $f_0$ and the $f_{k}(i,j)$'s. \ Denote 
$$
\begin{array}{l}
	A_{0}=\{c_0\} \\ \\ 
     A_{k}=\{c(k,l) \ | \ l=0,\ldots,p-1\} \ \ \ k=1,\ldots,h \\  
     A=\displaystyle{\bigcup_{k=0}^{h}A_k \ \ \ , \ \ \ A^{*}=\prod_{k=1}^{h}A_k} 
\end{array}
$$
\vskip 0.1cm\noindent
\textbf{1.10. Theorem.} \ $\mathit{S(\fg)^{\fg}=S_{p}(\fg)[A]\cong R/I}$.
\textit{Consequently,} $\mathit{S(\fg)^{\fg}}$ \textit{is a Cohen-Macaulay ring.}
\vskip 0.1cm
\textbf{Proof.} \ The ring $S_{p}(\fg)^{\fg}$ is integral over the ring $S_{p}(\fg)[A]$ (being so over $S_{p}(\fg)$).
Using (1), one has 
$$
	c(k,l)=c(k,1)^{l}/C(k)^{p(l-1)} \ \ \ l=2,\ldots,p-1, \ \ k=1\ldots h.
$$
Therefore $Q(S_{p}(\fg)[A])=Q(S_{p}(\fg)[c_{0},c(k,1) \ | \ k=1,\ldots,h])$.
By Theorem 1.7, $S_{p}(\fg)^{\fg}$ and  $S_{p}(\fg)[A]$ have the same quotient field.
It is therefore suffices to prove that the ring $S_{p}(\fg)[A]$ is normal.
\vskip 0.1cm
From relations (1) it follows that $S_{p}(\fg)[A]$ is generated by $A^{*}$ as a module over $S_{p}(\fg)[A_0]$, 
hence $A^{*}$ generates  $Q(S_{p}(\fg)[A])$ as a linear space over  $Q(S_{p}(\fg)[A_0])$.   
But $c(k,1)\not\in S_p(\fg)[A_0\cup\ldots\cup A_{k-1}]$ (see step 1 of the proof of 1.7) while $c(k,1)^{p}\in S_{p}(\fg)[A_0]$, $k=1,\ldots,h$, hence  $Q(S_{p}(\fg)[A])=Q(S_p(\fg)[c_{0},c(k,1) \ | \ k=1,\ldots,h])$ is of dimension $p^h$ over $Q(S_{p}(\fg)[A_0])$ with basis \\
$\{c(1,1)^{i_1}c(2,1)^{i_2}\cdots c(h,1)^{i_h} \ | \ i_j=0,\ldots,p-1\}$.
Thus the $p^h$-set $A^{*}$ form a basis of $Q(S_{p}(\fg)[A])$ as a linear space over  $Q(S_{p}(\fg)[A_0])$ and therefore form a free basis of $S_{p}(\fg)[A]$ over $S_{p}(\fg)[A_0]$.
In particular, $S_{p}(\fg)[A]$ is a Cohen-Macaulay ring (as a free module over the polynomial ring $S_{p}(\fg)[c_0]$, see 1.8), thereby satisfies the condition $(S_{2})$.   
\vskip 0.1cm
Since $c_0\not\in Q(S_{p}(\fg))$, the polynomial $f_0$ is a prime element of $S_{p}(\fg)[t_0]$, hence $S_{p}(\fg)[t_0]/(f_0)\cong S_p(\fg)[c_0]$.
Let $R_0=S_{p}(g)[t(k,i) \ | \ k=1,\ldots,h, \ i=1,\ldots,p-1]$ and let $I_0$ be the ideal of $R_0$ generated by the polynomials $f_{k}(i,j)$. 
From the definition of the $f_{k}(i,j)$'s (1.9) it follows that $R_0/I_0$ is generated by the homomorphic image of $\{1,t(k,i) \ | \ k=1,\ldots,h, \ i=1,\ldots,p-1\}$ as a module over $S_{p}(\fg)[A_0]$. 
Since $A^{*}$ is a free basis of $S_{p}(\fg)[A]$ over $S_{p}(\fg)[A_0]$, the natural map $R_0/I_0\to S_{p}[A]$ is a ring isomorphism. 
Hence $R/I\cong S_p[A]$.  
\vskip 0.1cm
We shall show that the ring $R/I$ satisfies the condition $(R_{1})$.
Firstly, from the isomorphism $R/I\cong S_p[A]$ we have $\Ht I=(p-1)h+1$.
Now, for $k=1,\ldots,h$ let 
$$
\begin{array}{l}
	g(k,1)=t(k,1)^{p}-c(k,1)^{p} \\  \\
     g(k,p-1)=t(k,p-1)^{p}-c(k,p-1)^{p}. 
\end{array}
$$
Then $g(k,1),g(k,p-1)\in I$. 
Indeed, it easy to verify that 
$$
\begin{array}{l}
	t(k,1)^{l}=t(k,l)C(k)^{(l-1)p} \ (\Mod I), \\  \\
     t(k,p-1)^{l}=t(k,p-l)M(k)^{(l-1)p} \ (\Mod I). 
\end{array}
$$
Thus 
$$
\begin{array}{ll}
	t(k,1)^{p}=t(k,1)t(k,1)^{p-1} & = t(k,0)M(k)^{p}C(k)^{(p-2)p} \\ \\
    									    & =(C(k)^{p-1}M(k))^{p}=c(k,1)^{p}	 \ (\Mod I),
\end{array}
$$
and
$$
\begin{array}{ll}	
      t(k,p-1)^{p}=t(k,p-1)t(k,p-1)^{p-1} & =t(k,0)M(k)^{p}M(k)^{(p-2)p} \\ \\ 
                                                              & =(C(k)M(k)^{p-1})^{p}=c(k,p-1)^{p} \ (\Mod I).
\end{array}
$$
Consider the following vectors: 
$$
\begin{array}{lll}
	\overrightarrow{f_{k}} & = & (f_{k}(1,1),f_{k}(1,2),\ldots,f_{k}(1,p-2))  \\ \\
     \overleftarrow{f_{k}} & = & (f_{k}(p-1,p-1),f_{k}(p-2,p-1),\ldots,f_{k}(2,p-1)) \\ \\ 
     \overrightarrow{t_{k}} & = & (t(k,2),t(k,3),\ldots,t(k,p-1)) \\ \\
     \overleftarrow{t_{k}} & = & (t(k,p-2),t(k,p-3),\ldots,t(k,1)) \\ \\
     \varphi & = & (f_0,\overrightarrow{f_{1}},g(1,1),\overrightarrow{f_{2}},g(2,1),\ldots,\overrightarrow{f_{h}},g(h,1)) \\ \\
     x & = & (e_{h+1,h+1}^{p},\overrightarrow{t_{1}},e_{1,1}^{p},\overrightarrow{t_{2}},e_{2,2}^{p},\ldots,\overrightarrow{t_{h}},e_{h,h}^{p}) \\ \\
     \varphi_{k} & = & \text{the vector which is obtained from} \ \varphi \ \text{by replacing} \ \overrightarrow{f_{k}} \ \text{with} \ \overleftarrow{f_{k}}, \\
     & & \text{and}\ g(k,1) \ \text{with} \ g(k,p-1). \\ \\
	x_k & = & \text{the vector which is obtained from} \ x \ \text{by replacing} \ \overrightarrow{t_{k}} \ \text{with} \ \overleftarrow{t_{k}}, \\
     & & \text{and}\ e_{k,k}^p \ \text{with} \ e_{k,n-k+1}^p.
\end{array}
$$
We have  
$$
\begin{array}{l}
\det(\partial{\varphi}/\partial{x})=-C(1)^{2p(p-1)}\cdots C(h)^{2p(p-1)} \\ \\
\det(\partial{\varphi_k}/\partial{x_k})=-C(k-1)^{p(2p-1)}M(k)^{2p(p-2)}T(k)^{p}+YC(k)^{p},
\end{array}
$$
where $Y\in S_{p}(\fg)$ and $C(0):=1$.
\\ \indent
Therefore, if $P$ is an element of the singular locus of $R$ such that $I\subset P$, then $P$ contains one of the $2$-sets $\{C(i)^{p},C(j)^{p}\}$, 
$\{C(k)^{p},T(k)^{p}\}$ (if $C(k)\in P$, consider $\det(\partial{\varphi_{k}}/\partial{x_{k}})$ and observe that $M(k)\in P$ if and only if $T(k)\in P$), hence contains a prime ideal of $S_{p}(\fg)$ of height $2$.  
Consequently, $P$ is of height $>\Ht{I}+1$. \eee
\vskip 0.1cm\noindent
\textbf{1.11.} \ Let $A$ be a commutative Noetherian ring (not necessarily local).
We say that $A$ is a \emph{complete intersection ring} if $A=R/I$ where $R$ is regular ring and $I$ is an ideal generated by an $R$-sequence.
This definition is customary in commutative algebra when $R$ is \emph{local}.
Our definition comes from algebraic geometry;
if $A$ is the coordinate ring of a projective variety over an algebraically closed field $k$,
then $A=R/I$ where $R$ is a polynomial ring over $k$ and the vanishing homogeneous ideal $I$ is of height $\codim{V}$.
Then $V$ is the intersection of $\codim{V}$ hypersurfaces and $I$ is generated by an $R$-sequence.
\\
\textbf{Corollary.} \ $\mathit{S(\fg)^{\fg}}$ \textit{is a complete intersection ring if and only if} $\mathit{p=2}$.
\vskip 0.1cm\noindent
\textbf{1.12. \ Question.} \ Is $S(\fg)^{\fg}$ a  Gorenstein ring? (for the definition of Gorenstein ring see e.g. [6]).
\vskip 0.1cm\noindent
\textbf{1.13.} \ Suppose that $p$ does not divide $n$.  
An explicit realization of the Poisson center $S(\BB)^{\BB}$ is next obtained from our realization of $S(\fg)^{\fg}$. \\
\indent
Clearly, $\fg=\BB\oplus Fc_0$.
Since $c_0$ is central, $S(\fg)^{\fg}=S(\BB)^{\BB}[c_0]$. 
The element $c_0$ is transcendental over $S(\BB)$, therefore, any generator of $S(\fg)^{\fg}$ has a unique expression as a polynomial in $c_0$ with coefficients in $S(\BB)^{\BB}$.
Furthermore, these coefficients generate $S(\BB)^{\BB}$ over $F$.   
\\ \indent
For $i=1,\ldots,n-1$, let $\varepsilon(i,i)=e_{i,i}-e_{n,n}$ be a basis element of the Cartan subalgebra of $\BB$. 
Clearly, the coefficients of $e_{i,j}^{p}$, $1\leq i\leq j\leq n$ (the generators of $S_{p}(\fg)$) belong to the polynomial ring
$$
	S_p(\BB)=F[\{\varepsilon(i,i)^{p} \ | \ 1\leq i\leq n-1\}\cup\{e_{i,j}^{p} \ | \ 1\leq i<j\leq n\}].
$$
\indent
Let us find the coefficients of $c(k,l)$, $1\leq k\leq h, \ 0\leq l\leq p-1$ (cf. 1.1, 1.2).
Since $D(k)\in\fg$, there exist $D_{\BB}(k)\in\BB$, $\alpha\in F$ such that $D(k)=D_{\BB}(k)+\alpha c_{0}$.
Explicitly (obtained by solving a system of linear equations), 
$$
D_{\BB}(k)=\left(1-\frac{2k-1}{n}\right)\sum_{i=1}^{k}\left(\varepsilon(i,i)+\varepsilon(n-i+1,n-i+1)\right)-\frac{2k-1}{n}\sum_{i=k+1}^{n-k}\varepsilon(i,i).
$$ 
Denote
$$
	M_{\BB}(k)=C(k)D_{\BB}(k)+T(k).
$$
From $c(k,l)=C(k)^{p-l}M(k)^{l}$ and $M(k)=C(k)D(k)+T(k)$ one get 
$$
	c(k,l)=C(k)^{p-l}(M_{\BB}(k)+\alpha C(k)c_0)^{l}.
$$ 
Denote 
$$
	c_{\BB}(k,l)=C(k)^{p-l}M_{\BB}(k)^{l}.
$$ 
By the binomial formula, 
$$
	c(k,l)=\sum_{i=0}^{l}\alpha^{l-i}\binom{l}{i}c_{\BB}(k,i)^{i}c_{0}^{l-i}.
$$ 
Therefore, up to scalar multiplication, the coefficients of the $c(k,l)$'s are the $c_{\BB}(k,l)$'s.  
\\ \indent
So we have
\vskip 0.1cm\noindent
\textbf{1.14. Theorem.} \ \textit{Suppose} $\mathit{p}$ \textit{does not divide} $\mathit{n.}$ 
\textit{The Poisson center} $\mathit{S(\BB)^{\BB}}$ \textit{is generated over} $\mathit{S_p(\BB)}$ \textit{by the elements} 
$\mathit{c_{\BB}(k,l)}$, $\mathit{1\leq k\leq h}$, $\mathit{1\leq l\leq p-1}$.
\vskip 0.1cm\noindent
\textbf{1.15. Theorem.} \ \textit{Suppose} $\mathit{p}$ \textit{does not divide} $\mathit{n}$.
\textit{The Poisson center} $\mathit{S(\BB)^{\BB}}$ \textit{is a Cohen-Macaulay ring.}
\vskip 0.1cm
\textbf{Proof.} \ $S(\fg)^{\fg}$ is a Cohen-Macaulay ring (Theorem 1.10) and  $S(\fg)^{\fg}=S(\BB)^{\BB}[c_0]$ where $c_0$ is transcendental over $S(\BB)^{\BB}$, 
so $S(\BB)^{\BB}$ is a Cohen-Macaulay ring.
\eee   
\subsection*{2. \ The center of $\mathbf{U(\fg)}$}
\textbf{2.1.} \ In this section, $\Char(F)=p>0$. Let $Z(\fg)$ be the center of the enveloping algebra $U(\fg)$.
According to 1.1, we use the same notation $e_{i,j}$ for a typical standard basis element of $\fg$, as well for $D(k)$,$C(k)$,$T(k)$ and $M(k)$, consider them as elements of $U(\fg)$.
The elements correspond to $c_0$ and $c(k,l)$ defined in 1.2 will respectively denoted by $z_0$ and $z(k,l)$.
\vskip 0.1cm\noindent
\textbf{2.2.} \ The corresponding polynomial ring of $S_p(\fg)$ in $Z(\fg)$ is
$$ 
	Z_{p}(\fg)=F[\{e_{i,i}^{p}-e_{i,i} \ | \ i=1,\ldots,n\}\cup\{e_{i,j}^{p} \ | \ 1\leq i<j\leq n\}].
$$
Clearly, $z_0\in Z(\fg)$.
The proof of Proposition 1.4, as is, shows that $C(k),M(k)\in U(\fg)^{\fn}$ and $z(k,l)\in Z(\fg)$. 
In Particular, $C(k)$ and $M(k)$ mutually commute and therefore relations (1) hold in $U(\fg)$, that is, 
\begin{equation} 
	z(k,i)z(k,j)=z(k,r(i,j))C(k)^{p(1-s(i,j))}M(k)^{ps(i,j)}.	
\end{equation}
\textbf{2.3.} \ The observations in 2.2 suggest that $Z(\fg)$ is generated over $Z_p(\fg)$ by $z_0$ and the $z(k,l)$'s, and it is isomorphic to $S(\fg)^{\fg}$ as a commutative algebra over $F$.
The following paragraphs is devoted to prove these facts.
\vskip 0.1cm\noindent
\textbf{2.4.} \ The well known facts in this paragraph are valid for \emph{any} finitely generated Lie algebra $\cL$ over \emph{any} field $F$ [2]. 
\\ \indent
For a non negative integer $m$, let $U_{m}(\cL)$ be the linear subspace of $U(\cL)$ generated by the products $x_{1}x_{2}\cdots x_{q}$, where $x_{1},\ldots,x_{q}\in\cL$ and $q\leq m$. 
If $u$ is a non zero element of $U(\cL)$, the smallest integer $m$ such that $u\in U_{m}(\cL)$ is termed the \emph{filtration of} $\mathit{u}$. 
If $u$ is a non zero element of $U(\cL)$ with filtration $m$, the \emph{grading of} $\mathit{u}$ is the image of $u$ in the linear space $U_{m}(\cL)/U_{m-1}(\cL)$ $(U_{-1}(\cL):=0)$ and we denoted it by $\gr{u}$. 
If $\mathcal{A}$ is a subalgebra of $U(\cL)$, the \emph{grading of} $\mathit{\mathcal{A}}$ is the commutative graded algebra
$$
	\gr{\mathcal{A}}=\bigoplus_{m\geq 0}(U_{m}(\cL)\cap\mathcal{A} + U_{m-1}(\cL))/U_{m-1}(\cL)
$$   
In particular, if $\mathcal{A}=U(\cL)$ then $\gr{U(\cL)}=\bigoplus_{m\geq 0}U_{m}(\cL)/U_{m-1}(\cL)\cong S(\cL)$.
\\
\textbf{Lemma.} \ If $u_1,\ldots,u_l\in\mathcal{A}$ such that $\gr{\mathcal{A}}=F[\gr{u_1},\ldots,\gr{u_l}]$, then $\mathcal{A}=F[u_1,\ldots,u_l]$.
\vskip 0.1cm\noindent
\textbf{2.5. Theorem.} \\
\textbf{a.} \ $\mathit{Z(\fg)=Z_{p}(\fg)[z_0,z(k,l) \ | \ 1\leq k\leq h, \ 1\leq l\leq p-1]}$. \\
\textbf{b.} \ $\mathit{Z(\fg)\cong S(\fg)^{\fg}}$.
\vskip 0.1cm
\textbf{Proof.} \ \textbf{a.} \ We set $S(\fg)=\gr{U(\fg)}$. Then $\gr{Z(\fg)}\subseteq S(\fg)^{\fg}$. 
Clearly, $\gr{u^{p}}=(\gr{u})^{p}$ for every element $u\in U(\fg)$, and $\gr(e_{i,i}^{p}-e_{i,i})=\gr{e_{i,i}^{p}}$ for each $i=1,\ldots n$.
By Theorem 1.10, $S(\fg)^{\fg}$ is generated over $F$ by the grading of the elements $e_{i,i}^{p}-e_{i,i}$ \ ($1\leq i\leq n$), \ $e_{i,j}^{p}$ \ ($1\leq i<j\leq n$), \ $z_0$ and $z(k,l)$ \ ($1\leq k\leq h$, $1\leq l\leq p-1$). Therefore,
$$
	S(\fg)^{\fg}=\gr(Z_{p}(\fg)[z_0,z(k,l) \ | \ 1\leq k\leq h, \ 1\leq l\leq p-1])\subseteq\gr{Z(\fg)}.
$$  
Hence $\gr{Z(\fg)}=S(\fg)^{\fg}$ and $\gr{Z(\fg)}$ is generated over $F$ by the grading of the following elements: \ $e_{i,i}^{p}-e_{i,i}$ \ ($1\leq i\leq n$), \ $e_{i,j}^{p}$ \ ($1\leq i<j\leq n$), \ $z_0$ and $z(k,l)$ \ ($1\leq k\leq h$, $1\leq l\leq p-1$). The assertion follows from Lemma 2.4. 
\\
\textbf{b.} \ By Theorem 1.10, $S(\fg)^{\fg}\cong R/I$. 
Let $\varphi:R\to Z(\fg)$ be the $F$-algebra epimorphism defined by $\varphi((\gr{e_{i,i}})^{p})=e_{i,i}^{p}-e_{i,i}$ \ ($1\leq i\leq n$), \ $\varphi(\gr(e_{i,j})^{p})=e_{i,j}^{p}$ \ ($1\leq i<j\leq n$), \
$\varphi(t_{0})=z_{0}$ and $\varphi(t(k,l))=z(k,l)$ \ ($1\leq k\leq h$, $1\leq l\leq p-1$). Then $I\subseteq\ker\varphi$ due to the relations in 2.2. Hence $Z(\fg)$ is a homomorphic image of $R/I$. 
The rings $Z(\fg)$ and $R/I$ are both domains with equal Krull dimension ($=K.\dim{S_{p}(\fg)}=K.\dim{Z_{p}(\fg)}=\dim_{F}\fg$), therefore $R/I\cong Z(\fg)$. \eee
\vskip 0.1cm\noindent
\textbf{2.6.} \ From Theorems 1.10, 2.5 we also have \\ 
\textbf{Corollary.} \ $\mathit{Z(\fg)}$ \textit{is a Cohen-Macaulay ring, and it is a complete intersection ring if and only if} $\mathit{p=2}$.
\vskip 0.1cm\noindent
\textbf{2.7.} \ In analogy to 1.12 we have \\
\textbf{Question.} \ Is $Z(\fg)$ a Gorenstein ring?.
\vskip 0.1cm\noindent
\textbf{2.8.} \ Suppose that $p$ does not divide $n$.  
An explicit realization of $Z(\BB)$ is next obtained from our realization of $Z(\fg)$. \\
\indent
Clearly, $\fg=\BB\oplus Fz_0$.
Since $z_0$ is central, $Z(\fg)=Z(\BB)[z_0]$. 
The element $z_0$ is transcendental over $Z(\BB)$, therefore, any generator of $Z(\fg)$ has a unique expression as a polynomial in $z_0$ with coefficients in $Z(\BB)$.
Furthermore, these coefficients generate $Z(\BB)$ over $F$.   
\\ \indent
For $i=1,\ldots,n-1$, let $\varepsilon(i,i)=e_{i,i}-e_{n,n}$ be a basis element of the Cartan subalgebra of $\BB$. 
The expression of a standard basis element of $\fg$ as a linear combination of $z_0$ and the basis elements of $\BB$, is over $\mathbb{F}_p$ (the prime field of $F$). 
It follows that the coefficients of the generators of $Z_p(\fg)$ belong to the polynomial ring
$$
	Z_p(\BB)=F[\{\varepsilon(i,i)^{p}-\varepsilon(i,i) \ | \ 1\leq i\leq n-1\}\cup\{e_{i,j}^{p} \ | \ 1\leq i<j\leq n\}].
$$
\indent
Let us find the coefficients of $z(k,l)$, $1\leq k\leq h, \ 0\leq l\leq p-1$.
Since $D(k)\in\fg$, there exist $D_{\BB}(k)\in \BB$, $\alpha\in F$ such that $D(k)=D_{\BB}(k)+\alpha z_{0}$.
Explicitly (obtained by solving a system of linear equations), 
$$
D_{\BB}(k)=\left(1-\frac{2k-1}{n}\right)\sum_{i=1}^{k}\left(\varepsilon(i,i)+\varepsilon(n-i+1,n-i+1)\right)-\frac{2k-1}{n}\sum_{i=k+1}^{n-k}\varepsilon(i,i).
$$ 
Denote
$$
	M_{\BB}(k)=C(k)D_{\BB}(k)+T(k).
$$
From $z(k,l)=C(k)^{p-l}M(k)^{l}$ and $M(k)=C(k)D(k)+T(k)$ one get 
$$
	z(k,l)=C(k)^{p-l}(M_{\BB}(k)+\alpha C(k)z_0)^{l}.
$$ 
Denote 
$$
	z_{\BB}(k,l)=C(k)^{p-l}M_{\BB}(k)^{l}.
$$ 
Since $C(k)$ and $M(k)$ mutually commute (Proposition 1.4a), the elements $C(k)$ and $M_{\BB}(k)$ mutually commute.
We can therefore apply the binomial formula, 
$$
	z(k,l)=\sum_{i=0}^{l}\alpha^{l-i}\binom{l}{i}z_{\BB}(k,i)^{i}z_{0}^{l-i}.
$$ 
Therefore, up to scalar multiplication, the coefficients of the $z(k,l)$'s are the $z_{\BB}(k,l)$'s.  
\\ \indent
So we have
\vskip 0.1cm\noindent
\textbf{2.9. Theorem.} \ \textit{Suppose} $\mathit{p}$ \textit{does not divide} $\mathit{n}$. 
\textit{The center} $\mathit{Z(\BB)}$ \textit{is generated over} $\mathit{Z_p(\BB)}$ \textit{by the elements} $\mathit{z_{\BB}(k,l)}$, $\mathit{1\leq k\leq h}$, $\mathit{1\leq l\leq p-1}$.
\vskip 0.1cm\noindent
\textbf{2.10. Theorem.} \ \textit{Suppose} $\mathit{p}$ \textit{does not divide} $\mathit{n}$. \textit{Then} $\mathit{S(\BB)^{\BB}\cong Z(\BB)}$.  
\textit{In particular}, $\mathit{Z(\BB)}$  \textit{is a Cohen-Macaulay ring.}
\vskip 0.1cm
\textbf{Proof.} \ Clearly, $S(\BB)^{\BB}\cong S(\BB)^{\BB}[c_0]/(c_0)$ and $Z(\BB)\cong Z(\BB)[z_0]/(z_0)$. 
Also, $S(\fg)^{\fg}=S(\BB)^{\BB}[c_0]$ and $Z(\fg)=Z(\BB)[z_0]$.
The isomorphism in Theorem 2.5b maps the variable $c_0$ to the variable $z_0$, hence induces an isomorphism of $S(\BB)^{\BB}$ onto $Z(\BB)$. 
By Theorem 1.15, $Z(\BB)$ is a Cohen-Macaulay ring.
\eee
\subsection*{3. \ From positive characteristic to zero characteristic}
\textbf{3.1.} \ Let $\cL$ be a finitely generated Lie algebra over a field $\KK$ of characteristic zero.
Suppose $\cL$ admits the following property:
\\
$\mathbf{P_1:}$ \ \textit{There exists a basis} $\mathit{B=(x_1,\ldots,x_r)}$ \textit{of} $\mathit{\cL}$ \textit{such that for each} $\mathit{i}$ \textit{the matrix} $\mathit{[\ad{x_{i}}]_{B}}$ \textit{of} $\mathit{\ad{x_i}}$
\textit{acting on} $\mathit{B}$ \textit{consists of integer entries.}
\\ \indent
Let $\phi$ be a field of a prime characteristic $p$.
Denote by $\cL_{p}(\phi)$ the "corresponding" algebra over $\phi$.
That is, $\cL_{p}(\phi)$ is the Lie algebra over $\phi$ with basis $C=(y_1,\ldots,y_r)$ such that $[\ad{y_{i}}]_{C}=[\ad{x_{i}}]_{B}~{(\!\!\!\!\mod{p})}$ \ for each $i$.
\\ \indent
For an element $a\in\ZZ[x_1,\ldots,x_r]$ denote by $\overline{a}$ the image of $a$ in $\FF_{p}[y_1,\ldots,y_r]$.
Thus $\overline{x_{i}}=y_{i}$ for each $i$.
\\ \indent
Let $\frak{h}$ be a Lie subalgebra of $\cL$ generated by a subset of $B$ and denote by  $\frak{h}_{p}(\phi)$ the corresponding Lie subalgebra over $\phi$.
Suppose we have following property as well: \\
$\mathbf{P_2:}$ \ \textit{There exist homogeneous polynomials} $\mathit{c_1,\ldots,c_s}$ \textit{in} $\mathit{S(\cL)^{\frak{h}}\cap\ZZ[x_1,\ldots,x_r]}$ \textit{such that from some prime} $\mathit{p}$ \textit{on, there exist homogeneous polynomials} $\mathit{q_1,\ldots,q_m\in\FF_{p}[y_1,\ldots,y_r]}$, $\mathit{\deg{q_i}\geq p}$ \textit{(with respect to the $y_i$'s) and} 
$$
    \mathit{S(\cL_{p}(\phi))^{\frak{h}_{p}(\phi)}=\phi[q_1,\ldots,q_m,\overline{c_1},\ldots,\overline{c_s}]}
$$
\textit{for every field} $\mathit{\phi}$ \textit{of characteristic} $\mathit{p}$.
\\ \indent
In property $\mathbf{P_2}$ it is crucial that $s$ and the invariants $c_{i}$ do not depend in $p$ (from some prime $p$ on).
Properties $\mathbf{P_1}$ and $\mathbf{P_2}$ suggest that $S(\cL)^{\frak{h}}=\KK[c_1,\ldots,c_s]$.
This is one of the main results of the present chapter.
\vskip 0.1cm
\noindent \textbf{3.2. Lemma. [1, Lemma 5.2]} \ \textit{Let} $\mathit{\mathcal{A}=\phi[a_1,\ldots,a_k]}$ \textit{be a subring of a polynomial ring} $\mathit{\phi[t_1,\ldots,t_n]}$ \textit{over a field} $\phi$
\textit{such that each} $\mathit{a_i}$ \textit{is homogenous of degree} $\mathit{d_{i}}$.
\textit{Assume} $\mathit{d_{i}\leq d_{i+1}}$ \textit{for each} $\mathit{i}$.
\textit{If} $\mathit{f\in\mathcal{A}}$ \textit{is homogenous of degree} $\mathit{d}$ \textit{(with respect to the} $\mathit{t_{i}}$\textit{'s) then} $\mathit{f\in\phi[a_1,\ldots,a_l]}$ \textit{where} $\mathit{d_l\leq d<d_{l+1}}$.
\vskip 0.1cm
\noindent\textbf{3.3.} \ For a $\KK$-subspace $M$ of $S(\cL)^{\frak{h}}$ we shall denote by $M_d$ the subspace of $M$ consisting of homogeneous elements of degree $d$. 
\\ \indent
The next theorem is a slight modification of [1, Theorem 5.3]. 
\\ 
\noindent\textbf{Theorem.} \  \textit{If} $\mathit{\cL}$ \textit{satisfies} $\mathit{\mathbf{P_1}}$ \textit{and} $\mathit{\mathbf{P_2}}$ \textit{then} $\mathit{S(\cL)^{\frak{h}}=\KK[c_1,\ldots,c_s]}$.
\vskip 0.1cm
\textbf{Proof.} \ Since $\cL$ satisfies $\mathbf{P_1}$, there exists a Lie algebra $\cL_{\QQ}$ having a Lie subalgebra  $\frak{h}_{\QQ}$ such that $\cL_{\QQ}\otimes_{\QQ}\KK=\cL$   and $\frak{h}_{\QQ}\otimes_{\QQ}\KK=\frak{h}$.
Since $S(\cL)=S(\cL_{\QQ})\otimes_{\QQ}\KK$  we have $S(\cL)^{\frak{h}}=S(\cL_{\QQ})^{\frak{h}_{\QQ}}\otimes_{\QQ}\KK$. 
Clearly, $S(\cL)^{\frak{h}}=\bigoplus_{d\geq 0}S(\cL)_{d}^{\frak{h}}$. 
Therefore, we can assume $\KK=\QQ$, $\cL=\cL_{\QQ}$,  $\frak{h}=\frak{h}_{\QQ}$ and we shall prove $S(\cL)^{\frak{h}}_{d}=\QQ[c_{1},\ldots,c_s]_d$ for every non negative integer $d$.
\\ \indent
Let $p$ be a prime number, $p>d$.
We have the natural epimorphism
$$
    \rho_{p}:\ZZ_p[x_1,\ldots,x_r]\longrightarrow\FF_p[y_1,\ldots,y_r]=S(\cL_{p}(\FF_p))
$$
considering $\ZZ_p[x_1,\ldots,x_r]$ as a subring of $S(\cL)=\QQ[x_1,\ldots,x_r]$ ($\ZZ_p$ is the localization of the ring of integers $\ZZ$ at $p$).
Note that $\rho_{p}$ is an extension of the natural map $\ZZ[x_1,\ldots,x_r]\to\FF_{p}[y_1,\ldots,y_r]$.
We clearly have
$$
\rho_p(S(\cL)^{\frak{h}}\cap \ZZ_p[x_1,\ldots,x_r])\subseteq S(\cL_{p}(\FF_p))^{\frak{h}_{p}(\FF_p)}.
$$
\indent Let $W$ be the $\ZZ_p$-submodule of $S(\cL)_{d}^{\frak{h}}$ consisting of all $\ZZ_p$-polynomials.
Let $V=\ZZ_p[c_1,\ldots,c_s]_{d}$.
Clearly, $V\subseteq W$ and they are finitely generated free $\ZZ_p$-modules.
\newline\indent
We shall next show that $V=W$.
Let $f\in W$.
Since $d<p$, it follows from Lemma 3.2 and $\mathbf{P_2}$ that $\rho_{p}(f)\in\FF_p[\overline{c_1},\ldots,\overline{c_{s}}]_{d}$.
Therefore $f\in\ZZ_p[c_1,\ldots,c_s]_d+p\ZZ_p[x_1,\ldots,x_r]$ and there exist $f_{0}\in V$, $h_1\in W$ such that $f=f_{0}+p h_1$.
Apply similar arguments to $h_1$ and conclude that $h_1=f_1+p h_2$ for some $f_1\in V$, $h_2\in W$.
Similarly we can successively choose $f_{i}\in V$ ($i=1,\ldots,n-1$) and $h_{n}\in W$ such that
$$
    f=f_{0}+p f_{1}+p^{2}f_{2}+\ldots+p^{n-1}f_{n-1}+p^{n}h_{n}.
$$
Denote $v_n=f_{0}+p f_{1}+p^{2}f_{2}+\ldots+p^{n-1}f_{n-1}$ for each $n$.
In particular $v_{n}\in V$.
Consider the completion $\widehat{W}$ with respect to the linear topology defined by the filtration $\{p^{n}W \ |n=1,2,\ldots \}$.
We identify $\widehat{V}$ as an $\ZZ_p$-submodule of $\widehat{W}$ (The topology of $\widehat{V}$ is the subspace topology;
is the same thing as the linear topology defined by the filtration $\{V\cap p^{n}W \ |n=1,2,\ldots \}$).
Let $\psi|W\to\widehat{W}$ be the natural map.
We have
$$
    \psi(f)=(f+p^{n}W)_{n=1}^{\infty}=(v_n+p^{n}h_{n}+p^{n}W)_{n=1}^{\infty}=
                (v_n+p^{n}W)_{n=1}^{\infty}\in\prod_{n=1}^{\infty}W/p^{n}W.
$$
Therefore, each component of $\psi(f)$ can be represented by $v_{n}\in V$.
Thus $\psi(f)\in\widehat{V}$.
Hence $\psi(W)\subseteq\widehat{V}$ because we have started with an arbitrary $f$ in $W$.
Hence $\widehat{W}\subseteq\widehat{V}$
(because $\widehat{W}=\Cl_{\widehat{W}}(\psi(W))\subseteq\Cl_{\widehat{W}}(\widehat{V})=\widehat{V}$, where $\Cl$ stands for the closure).
So we have $\widehat{V}=\widehat{W}$.
From the isomorphism $\widehat{W/V}\cong\widehat{W}/\widehat{V}$ ([6], Theorem 8.1) we get $\widehat{W/V}=0$.
Consider the $p$-adic completion $\widehat{\ZZ_p}$.
From the isomorphism $\widehat{W/V}\cong(W/V)\otimes_{\ZZ_p}\widehat{\ZZ_p}$ ([6], Theorem 7.2) we conclude $V=W$.
Finally, $S(\cL)_{d}^{\frak{h}}=W\otimes_{\ZZ_p}\QQ=V\otimes_{\ZZ_p}\QQ=\QQ[c_{1},\ldots,c_{r}]_{d}$.
\eee
\vskip 0.1cm
\noindent \textbf{3.4. Theorem.} \ \textit{Suppose that} $\mathit{p=\Char(F)=0}$. 
\textit{Then} $\mathit{S(\fg)^{\fg}=F[c_0]}$ \textit{and} $\mathit{Z(\fg)=F[z_0]}$.
\textit{In particular,} $\mathit{S(\fg)^{\fg}\cong Z(\fg)}$ \textit{as polynomial rings of one variable.}
\vskip 0.1cm
\textbf{Proof.} \ Clearly, $\fg$ satisfies $\mathbf{P_1}$ with respect to the standard basis $(e_{i,j})_{1\leq i\leq j\leq n}$.
By Theorem $1.10$, for every prime $p$ and for every field $\phi$ of characteristic $p$, the corresponding algebra $\fg_{p}(\phi)$ is generated over $\phi$ by the elements $c_0$, $e_{i,j}^{p}$ and $c(k,l)$, $1\leq k\leq\lfloor\frac{n-1}{2}\rfloor$, $1\leq l\leq p-1$. They are all homogeneous with respect to the standard basis and, the degrees of  $e_{i,j}^{p}$ and $c(k,l)$ are $\geq p$ $(\deg{c(k,l)}=kp+l)$.
Therefore, $\fg$ satisfies $\mathbf{P_2}$.
By Theorem 3.3 we have $S(\fg)^{\fg}=F[c_0]$. 
\\ \indent
We set $S(\fg)=\gr{U(\fg)}$ (see 2.4). Then $\gr{Z(\fg)}\subseteq S(\fg)^{\fg}$.
On the other hand, 
$$
	S(\fg)^{\fg}=F[c_0]=F[\gr{z_0}]=\gr{F[z_0]}\subseteq\gr{Z(\fg)}.
$$
Hence $\gr{Z(\fg)}=\gr{F[z_0]}$. By Lemma 2.4 we have $Z(\fg)=F[z_0]$.
\eee
\vskip 0.1cm
\noindent \textbf{3.5. Remark.} \ In Theorem 3.4, if  the field $F$ is \textit{algebraically closed}, the fact $Z(\fg)=F[z_0]$ also can be deduced from $S(\fg)^{\fg}=F[c_0]$ by using \textit{Duflo isomorphism} [4, Theorem 10.4.5]. 
\vskip 0.1cm
\noindent \textbf{3.6. Theorem.} \ \textit{Suppose that} $\mathit{p=\Char(F)=0}$. 
\textit{Then} $\mathit{S(\BB)^{\BB}=Z(\BB)=F}$.
\vskip 0.1cm
\textbf{Proof.} \ $S(\fg)=S(\BB)^{\BB}[c_0]$ and $Z(\fg)=Z(\BB)[c_0]$. 
By Theorem 3.4, $S(\BB)^{\BB}[c_0]=F[c_0]$ and $Z(\BB)[z_0]=F[z_0]$. 
Therefore $S(\BB)^{\BB}=Z(\BB)=F$. 
\eee
\vskip 0.1cm
\noindent \textbf{3.7. Remark.} \ Theorem 3.6 can also be deduced from Theorems 1.14 and 3.3; 
$\BB$ clearly admits $\mathbf{P_1}$ and $\mathbf{P_2}$ where the set of the $q_i$'s in 3.1 is the set of the $c_{\BB}(k,l)$'s, and the set of the $c_{i}'s$ in 3.1 is empty 
(loosely speaking, there are no zero characteristic generators while the degrees of the positive characteristic generators are getting higher as the characteristic grows). 
\subsection*{4. \ The Poisson semi-center of $\mathbf{S(\fg)}$}
\textbf{4.1.} \ Let $\cL$ be a finite dimensional Lie algebra over $F$.
A non zero element $f\in S(\cL)$ is called \textit{semi-invariant with weight} $\mathit{\lambda\in\cL^{*}}$ if $\ad{x}(f)=\lambda(x)f$ for all $x\in\cL$.
The $F$-algebra generated by the semi-invariants in $S(\cL)$ is denoted by $S(\cL)_{\si}^{\cL}$ and is called the \textit{Poisson semi-center of} $\mathit{S(\cL)}$.
The Poisson center $S(\cL)^{\cL}$ is therefore the $F$-subalgebra of $S(\cL)_{\si}^{\cL}$ consists of the semi-invariants with weight $0$.
\vskip 0.1cm
\noindent \textbf{4.2.} \ We shall use the following well known linear algebra result.
\vskip 0.1cm 
\noindent{\textbf{Lemma.}} \ Let $\frak{V}$ be a linear space over a field $\phi$, and $d_1,\ldots,d_r$ commuting linear transformations on $\frak{V}$.
Suppose that each $d_i$ satisfies a semi-simple split equation over $\phi$.
Let $\frak{U}$ be the linear subspace of $Hom_{\phi}(\frak{V})$ generated by $d_1,\ldots,d_r$.
For $\lambda\in\frak{U}^{*}$ let $\frak{V}_{\lambda}$ the linear subspace of $\frak{V}$ consisting of the elements $v$ such that $d_{i}v=\lambda(d_i)v$ for all $i=1,\ldots,r$.
Then $\frak{V}=\bigoplus{\frak{V}_{\lambda}}$ (a finite direct sum).
\vskip 0.1cm 
\noindent{\textbf{Proposition.}} \ If $p>0$ then $\mathit{S(\fg)_{\si}^{\fg}=S(\fg)^{\fn}}$.
\vskip 0.1cm
\textbf{Proof.} \ For a linear functional $\lambda\in\fg^{*}$, let $S(\fg)_{\lambda}=\{f\in S(\fg) \ | \ \ad{x}(f)=\lambda(x)f, \ x\in\fg\}$.
Then 
$$
	S(\fg)_{\si}^{\fg}=\bigoplus_{\lambda\in\fg^{*}}S(\fg)_{\lambda}.
$$
Let $\lambda\in\fg^{*}$.
Since $\fn=[\fg,\fg]$, $\lambda(\fn)=0$.
Therefore, $S(\fg)_{\lambda}\subseteq S(\fg)^{\fn}$, thus $S(\fg)_{\si}^{\fg}\subseteq S(\fg)^{\fn}$. 
\\ \indent
Apply the above Lemma with $\frak{V}=S(\fg)^{\fn}$ and $d_i=\ad{e_{i,i}}$ and conclude $S(\fg)^{\fn}\subseteq S(\fg)_{\si}^{\fg}$ (the transformations $\ad{e_{i,i}}$ satisfy the polynomial $t^{p}-t$).
\eee
\vskip 0.1cm
\noindent \textbf{4.3. Theorem.} \ \textit{If} $\mathit{p>0}$ \textit{then} 
$$
	\mathit{Q(S(\fg)_{\si}^{\fg})=Q(S_{p}(\fg)[c_0,C(k),M(k) \ | \ k=1,\ldots,\lfloor\frac{n}{2}\rfloor])}.
$$ 
\indent\textbf{Proof.} \ Clearly, the field extension $Q(S_{p}(\fg)[c_0])/Q(S_p(\fg))$ is of degree $p$. \\ 
The extension $Q(S_{p}(\fg)[c_0,C(1)])/Q(S_p(\fg)[c_0])$ is of degree $p$ because $\ad{e_{1,1}}(C(1))=C(1)$ while $\ad{e_{1,1}}(S_{p}(\fg)[c_0])=0$.  
The extension $Q(S_{p}(\fg)[c_0,C(1),M(1)])/Q(S_p(\fg)[c_0,C(1)])$ is of degree $p$ since $\ad{e_{n,1}}(M(1))\not\in S(\fg)$ while $\ad{e_{n,1}}(S_p(\fg)[c_0,C(1)])\in S(\fg)$ (as in the proof of Theorem 1.7). 
\\ \indent
For $1\leq k<\lfloor\frac{n}{2}\rfloor$, the extension 
$$
	Q(S_p(\fg)[c_0,C(l),M(l),C(k+1) \ | \ l=1,\ldots,k])/Q(S_p(\fg)[c_0,C(l),M(l) \ | \ l=1,\ldots,k])
$$
is of degree $p$ since $\ad{e_{k+1,k+1}}(C(k+1))=C(k+1)$ while \\ $\ad{e_{k+1,k+1}}(S_p[c_0])=0$, $\ad{e_{k+1,k+1}}(C(l))=0$ and $\ad{e_{k+1,k+1}}(M(l))=0$ for $1\leq l\leq k$.
\\ \indent
If $n$ is even, the extension 
$$
	Q(S_p(\fg)[c_0,C(l),M(l) \ | \ l=1,\ldots,k+1])/Q(S_p(\fg)[c_0,C(l),M(l),C(k+1) \ | \ l=1,\ldots,k])	
$$
is of degree $p$ because $\ad{e_{n-k,k+1}}(M(k+1))\not\in S(\fg)$ while $\ad{e_{n-k,k+1}}(C(k+1))\in S(\fg)$, \\
$\ad{e_{n-k,k+1}}(S_p(\fg)[c_0])=0$, $\ad{e_{n-k,k+1}}(C(l))=0$ and $\ad{e_{n-k,k+1}}(M(l))=0$ for $1\leq l\leq k$.
It follows that the field extension 
$$
	Q(S_{p}(\fg)[c_0,C(k),M(k) \ | \ k=1,\ldots,\lfloor\frac{n}{2}\rfloor])/Q(S_p(\fg))
$$ 
is of degree $p^{1+\lfloor\frac{n}{2}\rfloor+h}=p^{n}$.
\\ \indent
We proceed as in [1, Proposition 1.13] and Theorem 1.7.
Let $r,i,j,k$ be integers such that $1\leq r\leq n-1$, $1\leq k\leq\lfloor\frac{n}{2}\rfloor$, $1\leq i\leq k$, $k+1\leq j\leq n-k$. 
For each $e_{i,j}\in \hat{d}(r)$ we have \\
(1) \ If $e_{i,j}$ is lying above the anti-diagonal of $\fm$ then 
$$
\begin{array}{l}
\ad{T_{i}(i,j)}(e_{i,j})\neq 0 \ \text{and} \\ \\ 
\ad{T_{i}(i,j)}(\hat{d}(r)\setminus\{e_{i,j}\}\cup \hat{d}(r+1)\cup\ldots\cup \hat{d}(n-1))=0.
\end{array} 
$$
(2) \  If $e_{i,j}$ is lying below the anti-diagonal of $\fm$, then 
$$
\begin{array}{l}
\ad{S_{i}(1,n-j+1)}(e_{i,j})\neq 0 \  \text{and} \\ \\
\ad{S_{i}(1,n-j+1)}(\hat{d}(r)\setminus\{e_{i,j}\}\cup \hat{d}(r+1)\cup\ldots\cup \hat{d}(n-1))=0. 
\end{array}
$$
Let $\{\fhn\}=\bigcup_{r=1}^{n-1}\hat{d}(r)$.
The set $\{\fhn\}$ consisting of $\frac{n(n-1)}{2}-\lfloor\frac{n}{2}\rfloor$ elements. 
It follows that the field extension $Q(S(\fg)^{\fg}[\{\fhn\}])/Q(S(\fg)^{\fg})$ is of degree $p^{\frac{n(n-1)}{2}-\lfloor\frac{n}{2}\rfloor}$. 
\\ \indent
Now, for $k=1,\ldots,\lfloor\frac{n}{2}\rfloor$ we have $\ad{C(k)}(e_{k,k})=-C(k)$, $\ad{C(k)}(e_{l,l})=0$ when $l<k$, and by Proposition 1.2.a. $\ad{C(k)}(S(\fg)^{\fn}[\{\fhn\}])=0$.
It follows that the field extension 
$$
	Q(S(\fg)^{\fn}[\{\fhn\}\cup\{e_{k,k} \ | \ k=1,\ldots,\lfloor\frac{n}{2}\rfloor\}])/Q(S(\fg)^{\fn}[\{\fhn\}])
$$ is of degree $p^{\lfloor\frac{n}{2}\rfloor}$.
Therefore, the degree of 
$$
	Q(S(\fg)^{\fn}[\{\fhn\}\cup\{e_{k,k} \ | \ k=1,\ldots,\lfloor\frac{n}{2}\rfloor\}])/Q(S_p(\fg))
$$ is at least 
$$
	p^{n}p^{\frac{n(n-1)}{2}-\lfloor\frac{n}{2}\rfloor}p^{\lfloor\frac{n}{2}\rfloor}=p^{\frac{n(n+1)}{2}}=p^{\dim{\fg}},
$$
which is the degree of $Q(S(\fg))/Q(S_{p}(\fg))$. 
By degree considerations we must have 
$$
	Q(S(\fg)^{\fn})=Q(S_{p}(\fg)[c_0,C(k),M(k) \ | \ k=1,\ldots,\lfloor\frac{n}{2}\rfloor]).
$$ 
Proposition 4.2 complete the proof. 
\eee
\vskip 0.1cm
\noindent \textbf{4.4.} \ Let $t_0, t_{C(k)},t_{M(l)}$,\ $k=1,\ldots,\lfloor\frac{n}{2}\rfloor$, $l=1,\ldots,h$, be algebraically independent elements over  $S_{p}(\fg)$.
Consider the polynomial ring (of $n$ variables):
$$
	R=S_{p}(g)[t_0,t_{C(k)},t_{M(l)} \ | \ k=1,\ldots\lfloor\frac{n}{2}\rfloor, \ l=1,\ldots h] 
$$
and the following $n$ polynomials of $R$: 
$$
	f_0=t_{0}^{p}-c_{0}^{p} \ \ , \ \ f_{C(k)}=t_{C(k)}^{p}-C(k)^{p} \ \ , \ \ f_{M(l)}=t_{M(l)}^{p}-M(l)^{p}. 
$$
Let $I$ be the ideal of $R$ generated by $f_0$, $f_{C(k)}$, $f_{M(l)}$, $k=1,\ldots,\lfloor\frac{n}{2}\rfloor$, $l=1,\ldots,h$. 
\vskip 0.1cm
\noindent \textbf{4.5. Theorem.} \ \textit{If} $\mathit{p>0}$ \textit{then} 
$$
	\mathit{S(\fg)_{\si}^{\fg}=S_{p}(\fg)[c_0,C(k),M(k) \ | \ k=1,\ldots,\lfloor\frac{n}{2}\rfloor]\cong R/I}.
$$ 
\textit{Consequently,} $\mathit{S(\fg)_{\si}^{\fg}}$ \textit{is a complete intersection ring.}
\vskip 0.1cm
\textbf{Proof.} \ Besides the use of the Jacobian criterion for regularity, our proof consists of identical arguments appear in [1, Theorem 1.21] (some of them also appear in the proof of Theorem 1.10] and we shall not repeat them. 
\\ \indent
Consider the following vectors (of length $n$):
$$
\varphi=\left\{\begin{array}{ll}
				(f_0,f_{C(1)},f_{M(1)},\ldots,f_{C(h)},f_{M(h)}) & n \ \text{is odd} \\ \\
				(f_0,f_{C(1)},f_{M(1)},\ldots,f_{C(h)},f_{M(h)},f_{C(h+1)}) & n \ \text{is even}
\end{array}\right.
$$
$$
x=\left\{\begin{array}{ll}
			(e_{h+1,h+1}^{p},e_{1,n}^{p},e_{1,1}^{p},e_{2,n-1}^{p},\ldots,e_{k,k}^{p},e_{k+1,n-k}^{p},\ldots,e_{h-1,h-1}^{p},e_{h,n-h+1}^{p},e_{h,h}^{p}) & n \ \text{is odd} \\ \\	
			(e_{h+1,h+1}^{p},e_{1,n}^{p},e_{1,1}^{p},e_{2,n-1}^{p},\ldots,e_{k,k}^{p},e_{k+1,n-k}^{p},\ldots,e_{h,h}^{p},e_{h+1,n-h}^{p}) & n \ \text{is even} 			
\end{array}\right.
$$
Then
$$
	\det(\partial{\varphi}/\partial{x})=\pm C(1)^{2p}\cdots C(h-1)^{2p}C(h)^{ip}	
$$
where $i=1$ if $n$ is odd, while $i=2$ if $n$ is even.
\\ \indent
Let $x_{k}$ be the vector which is obtained from $x$ by replacing $e_{k,k}^{p}$ with $e_{k,k+1}^{p}$, and $e_{k+1,n-k}^{p}$ with $e_{k,n-k}^{p}$, $1\leq k\leq h$ (if $k=h$ and $n$ is odd, we only replace $e_{h,h}^{p}$ with $e_{h,h+1}^{p}$).
We have
$$
	\det(\partial{\varphi}/\partial{x_k})=\pm C(1)^{2p}\cdots C(k-1)^{2p}T_{k}(k,k+1)^{\alpha p}C(k+1)^{2p}\cdots C(h-1)^{2p}C(h)^{\beta p}
$$
where $\alpha,\beta\in\{1,2\}$.
Therefore, if $P$ is an element of the singular locus of $R$ such that $I\subset P$, then $P$ contains one of the $2$-sets $\{C(i)^{p},C(j)^{p}\}$, 
$\{C(k)^{p},T_{k}(k,k+1)^{p}\}$, hence contains a prime ideal of $S_{p}(\fg)$ of height $2$.
\eee
\vskip 0.1cm
\noindent \textbf{4.6. Theorem.} \ \textit{Suppose that} $\mathit{\Char(F)=p=0}$.
\textit{Then} 
$$
	\mathit{S(\fg)_{\si}^{\fg}=F[c_0,C(k),M(k) \ | \ k=1,\ldots,\lfloor\frac{n}{2}\rfloor]}.
$$
\textit{In particular, the Poisson center} $\mathit{S(\fg)_{\si}^{\fg}}$ \textit{is a polynomial ring of} $\mathit{n}$ \textit{variables.}
\vskip 0.1cm
\textbf{Proof.} \ The assertions follows from 3.3, 4.2, 4.5 and PBW.
\eee 
\vskip 0.1cm
\noindent \textbf{4.7.} \ Suppose that $p$ does not divide $n$.  
An explicit realization of the Poisson semi-center $S(\BB)_{\si}^{\BB}$ is next obtained from our realization of $S(\fg)_{\si}^{\fg}$. \\
\indent
Clearly, $\fg=\BB\oplus Fc_0$ and $c_0$ is central.
By extending $\lambda\in\BB^{*}$ to $\fg^{*}$ by $\lambda(c_0)=0$ one get $S(\BB)_{\si}^{\BB}[c_0]\subseteq S(\fg)_{\si}^{\fg}$.
Since $c_0$ is transcendental over $S(\BB)$ we have $S(\fg)_{\si}^{\fg}\subseteq S(\BB)_{\si}^{\BB}[c_0]$. 
Hence $S(\fg)_{\si}^{\fg}=S(\BB)_{\si}^{\BB}[c_0]$ (thus $S(\BB)_{\si}^{\BB}=S(\BB)^{\fn}$) 
and any generator of $S(\fg)_{\si}^{\fg}$ has a unique expression as a polynomial in $c_0$ with coefficients in $S(\BB)_{\si}^{\BB}$.
Furthermore, these coefficients generate $S(\BB)_{\si}^{\BB}$ over $F$. 
By Theorem 4.5 (and with the notation of 1.13) we therefore have
\vskip 0.1cm
\noindent \textbf{4.8. Theorem.} \  \textit{Suppose} $\mathit{p>0}$ \textit{and} $\mathit{p}$ \textit{does not divide} $\mathit{n}$. \textit{Then} 
$$
	\mathit{S(\BB)_{\si}^{\BB}=S_{p}(\BB)[C(k),M_{\BB}(k) \ | \ k=1,\ldots,\lfloor\frac{n}{2}\rfloor]},
$$ 
\textit{and} $\mathit{S(\BB)_{\si}^{\BB}}$ \textit{is a complete intersection ring.}
\vskip 0.1cm 
\noindent \textbf{4.9.} \ By Theorems 3.3, 4.8 we have 
\\
\textbf{Theorem.} \ \textit{Suppose that} $\mathit{\Char(F)=p=0}$.
\textit{Then} 
$$
	\mathit{S(\BB)_{\si}^{\BB}=F[C(k),M_{\BB}(k) \ | \ k=1,\ldots,\lfloor\frac{n}{2}\rfloor]}.
$$
\textit{In particular, the Poisson center} $\mathit{S(\BB)_{\si}^{\BB}}$ \textit{is a polynomial ring of} $\mathit{n-1}$ \textit{variables.}
\subsection*{5. \ The semi-center of $\mathbf{U(\fg)}$}
\textbf{5.1.} \ \ Let $\cL$ be a finite dimensional Lie algebra over $F$.
A non zero element $u\in Z(\cL)$ is called \textit{semi-central with weight} $\mathit{\lambda\in\cL^{*}}$ if $[x,u]=\lambda(x)u$ for all $x\in\cL$.
The $F$-algebra generated by the semi-central elements in $Z(\cL)$ is denoted by $Sz(\cL)$ and is called the \textit{semi-center of} $\mathit{U(\cL)}$.
Clearly, $Z(\cL)$ is  the $F$-subalgebra of $Sz(\cL)$ consists of the semi-central elements with weight $0$.
By [3, Proposition 2.1] $Sz(\cL)$ is commutative.
\vskip 0.1cm
\noindent \textbf{5.2.} \ We use the same notation for $C(k)$,$M(k)$ and $M_{\BB}(k)$, consider them as elements of $U(\fg)$.
\vskip 0.1cm
\noindent\textbf{5.3.} \ The arguments here follow those of Theorems 2.5, 3.4 almost verbatim, with $C(k)$, $M(k)$ replacing $c(l,k)$.
\vskip 0.1cm\noindent 
\textbf{Theorem.} \\
\textbf{a.} \ \textit{If} $\mathit{p>0}$ \textit{then} 
$$
	\mathit{Sz(\fg)=Z_{p}(\fg)[z_0,C(k),M(k) \ | \ k=1,\ldots,\lfloor\frac{n}{2}\rfloor]\cong S(\fg)_{\si}^{\fg}}.
$$ 
\indent\textit{In particular,} $\mathit{Sz(\fg)}$ \textit{is a complete intersection ring.} \\
\textbf{b.} \ \textit{If} $\mathit{p=0}$ \textit{then}
$$
	\mathit{Sz(\fg)=F[z_0,C(k),M(k) \ | \ k=1,\ldots,\lfloor\frac{n}{2}\rfloor]\cong S(\fg)_{\si}^{\fg}}.
$$ 
\indent\textit{In particular,} $\mathit{Sz(\fg)}$ \textit{is a polynomial ring of $n$ variables.}
\vskip 0.1cm
\noindent\textbf{5.4.} \ The arguments here follow those of Theorems 2.9, 2.10, 3.6 almost verbatim, with $C(k)$, $M_{\BB}(k)$ replacing $c_{\BB}(l,k)$.
\vskip 0.1cm\noindent 
\textbf{Theorem.} \\
\textbf{a.} \ \textit{Suppose} $\mathit{p>0}$ \textit{and} $\mathit{p}$ \textit{does not divide} $\mathit{n}$. \textit{Then} 
$$
	\mathit{Sz(\BB)=Z_{p}(\BB)[C(k),M_{\BB}(k) \ | \ k=1,\ldots,\lfloor\frac{n}{2}\rfloor]\cong S(\BB)_{\si}^{\BB}}.
$$ 
\indent\textit{In particular,} $\mathit{Sz(\BB)}$ \textit{is a complete intersection ring.} \\
\textbf{b.} \ \textit{If} $\mathit{p=0}$ \textit{then}
$$
	\mathit{Sz(\BB)=F[C(k),M_{\BB}(k) \ | \ k=1,\ldots,\lfloor\frac{n}{2}\rfloor]\cong S(\BB)_{\si}^{\BB}}.
$$ 
\indent\textit{In particular,} $\mathit{Sz(\BB)}$ \textit{is a polynomial ring of $n-1$ variables.}
\subsection*{Acknowledgments}
\ \ \ \ \ The author wishes to thank to his mentor, professor Amiram Braun, for extensive help in the preparation of this paper.  

\end{document}